# An arguable addition to the standard Deduction Theorems of first order theories[1]

Bhupinder Singh Anand[2]

We consider an arguable addition to the standard Deduction Theorems of first order theories.

## Contents[3]

---

[1] *Updated: Monday 13$^{th}$ September 2004 3:31:36 AM IST by re@alixcomsi.com*

[2] The author is an independent scholar. E-mail: re@alixcomsi.com; anandb@vsnl.com.

[3] Key words: Arithmetic, Church, classical, computable, constructive, deduction, effective, expressible, finite, formal, formula, function, Gödel, interpretation, meta-assertion, meta-theorem, natural number, number-theoretic, numeral, Peano, primitive, proof, recursive, relation, satisfiable, sequence, standard, truth, undecidable.

## 1. Introduction

We first review, in Meta-theorem 1, the proof of a standard Deduction Theorem - ([$T$], [$A$])|-$_K$ [$B$] if, and only if, [$T$]|-$_K$ [$A => B$] (cf. [Me64], Corollary 2.6, p61) - of classical[4], first-order, theories, where an explicit deduction of [$B$][5] from ([$T$], [$A$]) is known.

We then show, in Corollary 1.2, that, assuming Church's Thesis, Meta-theorem 1 can be constructively extended to cases where, ([$T$], [$A$])|-$_K$ [$B$] is established meta-mathematically, assuming the consistency of ([$T$], [$A$]), but where an explicit deduction of [$B$] from ([$T$], [$A$]) is not known.

We finally argue, in Meta-theorem 2, that:

   ([$T$], [$A$])|-$_K$ [$B$] holds if, and only if, [$T$]|-$_K$ [$B$] holds when we assume [$T$]|-$_K$ [$A$].

(In other words, that [$B$] is a deduction from ([$T$], [$A$]) in K if, and only if, whenever [$A$][6] is a hypothetical deduction from [$T$] in K, [$B$] is a deduction from [$T$] in K.)

## 2. A standard Deduction Theorem

The following is, essentially, Mendelson's proof of a standard Deduction Theorem ([Me64], p61, Proposition 2.4 and Corollary 2.6) of an arbitrary first order theory K:

---

[4] We take Mendelson [Me64] as representative, in the area that it covers, of standard expositions of classical, first order, logic.

[5] We use square brackets to differentiate between a formal expression [$F$] and its interpretation "$F$", where we follow Mendelson's definition of an interpretation M of a formal theory K, and of the interpretation of a formula of K under M ([Me64], p49, §2). For instance, we use [$n$] to denote the numeral in K whose standard interpretation is the natural number $n$.

[6] For the purposes of this essay, we assume everywhere that [$T$] is an abbreviation for a finite set of K-formulas {[$T_1$], [$T_2$], ... [$T_l$]}, whereas [$A$], [$B$], ... are closed well-formed formulas of K. We note, also, that "[$A$] & [$B$]" and "[$A$ & $B$]" denote the same K-formula.

**Meta-theorem 1**: If $[T]$ is a set of well-formed formulas of an arbitrary first order theory K, and if $[A]$ is a closed well-formed formula of K, and if $([T], [A])|\text{-}_K [B]$, then $[T]|\text{-}_K [A => B]$[7].

**Proof**: Let $<[B_1], [B_2], ..., [B_n]>$ be a deduction of $[B]$ from $([T], [A])$ in K.

Then, by definition, $[B_n]$ is $[B]$ and, for each $i$, either $[B_i]$ is an axiom of K, or $[B_i]$ is in $[T]$, or $[B_i]$ is $[A]$, or $[B_i]$ is a direct consequence by some rules of inference of K of some of the preceding well-formed formulas in the sequence.

We now show, by induction, that $[T]|\text{-}_K [A => B_i]$ for each $i =< n$. As inductive hypothesis, we assume that the proposition is true for all deductions of length less than $n$.

(i) If $[B_i]$ is an axiom, or belongs to $[T]$, then $[T]|\text{-}_K [A => B_i]$, since $[B_i => (A => B_i)]$ is an axiom of K.

(ii) If $[B_i]$ is $[A]$, then $[T]|\text{-}_K [A => B_i]$, since $[T]|\text{-}_K [A => A]$.

(iii) If there exist $j$, $k$ less than $i$ such that $[B_k]$ is $[B_j => B_i]$, then, by the inductive hypothesis, $[T]|\text{-}_K [A => B_j]$, and $[T]|\text{-}_K [A => (B_j => B_i)]$. Hence, $[T]|\text{-}_K [A => B_i]$.

(iv) Finally, suppose there is some $j < i$ such that $[B_i]$ is $[(Ax)B_i]$, where $x$ is a variable in K. By hypothesis, $[T]|\text{-}_K [A => B_j]$. Since $x$ is not a free variable of $[A]$, we have that $[(Ax)(A => B_j) => (A => (Ax)B_j]$ is PA-provable. Since $[T]|\text{-}_K [A => B_j]$, it follows by Generalisation that $[T]|\text{-}_K [(Ax)(A => B_j)]$, and so $[T]|\text{-}_K [A => (Ax)B_j]$, i.e. $[T]|\text{-}_K [A => B_i]$.

---

[7] The converse is trivially true (cf. [Sh67], p33).

This completes the induction, and Meta-theorem 1 follows as the special case where $i = n$. ¶[8]

### 2.1 A number-theoretic corollary

Now, Gödel has defined ([Go31], p22, Definition 45(6)) a primitive recursive number-theoretic relation $xB_{(K, [T])}y$ that holds if, and only if, $x$ is the Gödel-number of a deduction from $T$ of the K-formula whose Gödel-number is $y$.

We thus have:

**Corollary 1.1**[9]: If the Gödel-number of the well-formed K-formula $[B]$ is $b$, and that of the well-formed K-formula $[A => B]$ is $c$, then Meta-theorem 1 holds if, and only if[10]:

$(Ex)xB_{(K, [T], [A])}b => (Ez)zB_{(K, [T])}c$

### 2.2 An extended Deduction Theorem

We next consider the meta-proposition:

**Corollary 1.2**: If we assume Church's Thesis[11], then Meta-theorem 1 holds even if the premise $([T], [A])|-_K [B]$ is established meta-mathematically, assuming the consistency of $([T], [A])$, but a deduction $<[B_1], [B_2], ..., [B_n]>$ of $[B]$ from $([T], [A])$ in K is not known explicitly.

---

[8] We use the symbol "¶" as an end-of-proof marker.

[9] We note that Corollary 1.1 and Corollary 2.2 may be essentially different number-theoretic assertions, which, in the absence of a formal proof, cannot be assumed to be equivalent.

[10] We note that this symbolically expresses a meta-equivalence in a recursive arithmetic RA, based on a semantic interpretation of the definition of the primitive recursive relation $xB_{(K, [T])}y$; it is not a K-formula.

[11] Church's Thesis: A number-theoretic function is effectively computable if, and only if, it is recursive ([Me64], p147, footnote). We appeal explicitly to Church's Thesis here to avoid implicitly assuming that every recursive relation is algorithmically decidable.

**Proof**: Since Gödel's number-theoretic relation $xB_{(K, [T])}y$ is primitive recursive, it follows that, if we assume Church's Thesis - which implies that a number-theoretic relation is decidable if, and only if, it is recursive - we can effectively determine some finite natural number $n$ for which the assertion $nB_{(K, [T], [A])}b$ holds, where the Gödel-number of the well-formed K-formula $[B]$ is $b$.

Since $n$ would then, by definition, be the Gödel-number of a deduction $<[B_1], [B_2], ..., [B_n]>$ of $[B]$ from $([T], [A])$ in K, we may thus constructively conclude, from the meta-mathematically determined assertion $([T], [A])|-_K [B]$, that some deduction $<[B_1], [B_2], ..., [B_n]>$ of $[B]$ from $([T], [A])$ in K can, indeed, be effectively determined. Meta-theorem 1 follows. ¶

## 3. An additional Deduction Theorem

We, finally, argue that:

**Meta-theorem 2**: If K is an arbitrary first order theory, and if $[A]$ is a closed well-formed formula of K, then $([T], [A])|-_K [B]$ holds if, and only if, $[T]|-_K [B]$ holds when we assume that $[T]|-_K [A]$ holds.

**Proof**: First, if there is a deduction $<[B_1], [B_2], ..., [B_n]>$ of $[B]$ from $([T], [A])$ in K, and there is a deduction, $<[A_1], [A_2], ..., [A_m]>$, of $[A]$ from $[T]$ in K, then $<[A_1], [A_2], ..., [A_m], [B_1], [B_2], ..., [B_n]>$ is a deduction of $[B]$ from $[T]$ in K. Hence, if $([T], [A])|-_K [B]$ holds, then $[T]|-_K [B]$ holds when we assume $[T]|-_K [A]$.

Second, if there is a deduction $<[B_1], [B_2], ..., [B_n]>$ of $[B]$ from $[T]$ in K, then we have, trivially, that, if $[T]|-_K [B]$ holds when we assume $[T]|-_K [A]$, then $([T], [A])|-_K [B]$ holds.

Last, we assume that there is no deduction, $<[B_1], [B_2], ..., [B_n]>$, of $[B]$ from $[T]$ in K. If, now, $[T]|-_{K'} [B]$ holds when we assume that $[T]|-_{K'} [A]$ holds in any extension K' of K,

then, if we assume that there is a sequence $<[A_1], [A_2], ..., [A_m]>$ of well-formed K'-formulas such that $[A_m]$ is $[A]$ and, for each $m >= i >= 1$, either $[A_i]$ is an axiom of K', or $[A_i]$ is in $[T]$, or $[A_i]$ is a direct consequence by some rules of inference of K' of some of the preceding well-formed formulas in the sequence, then it follows from our hypothesis[12] that we can show, by induction on the deduction length $n$, that there is a sequence $<[B_1], [B_2], ..., [B_n]>$ of well-formed K-formulas such that $[B_1]$ is $[A]$[13], $[B_n]$ is $[B]$ and, for each $i > 1$, either $[B_i]$ is an axiom of K, or $[B_i]$ is in $[T]$, or $[B_i]$ is a direct consequence by some rules of inference of K of some of the preceding well-formed formulas in the sequence.

In other words, if there is a deduction, $<[A_1], [A_2], ..., [A_m]>$, of $[A]$ from $[T]$ in K', then, by our hypothesis, $<[A], [B_2], ..., [B_n]>$ is a deduction of $[B]$ from $([T], [A])$ in K'. By definition, it follows that $<[A], [B_2], ..., [B_n]>$ is, then, a deduction of $[B]$ from $([T], [A])$ in K. We thus have that, if $[T]|-_K [B]$ holds when we assume $[T]|-_K [A]$, then $([T], [A])|-_K [B]$ holds. This completes the proof. ¶

In view of Corollary 1.2, we thus have:

**Corollary 2.1**: If we assume Church's Thesis, and if $[A]$ is a closed well-formed formula of K, then we may conclude that $[T]|-_K ([A] => [B])$ holds if $[T]|-_K [B]$ holds when we assume $[T]|-_K [A]$.[14]

We note that, in the notation of Corollary 1.1, if the Gödel-number of the well-formed K-formula $[A]$ is $a$, then Corollary 2.1 holds if, and only if[15]:

---

[12] That $[T]|-_K [A]$ holds.

[13] $[A]$ is, thus, the hypothesis in the sequence; it is the only well-formed K-formula in the sequence that is not an axiom of K, not in $[T]$, and not a direct consequence of the axioms of K by any rules of inference of K.

[14] We give a model-theoretic proof of Corollary 2.1 in the Appendix.

[15] We note that this, too, is not a K-formula, but a semantic meta-equivalence, based on the definition of the primitive recursive relation $xB_{(K, [T])}y$.

**Corollary 2.2**: $((Ex)xB_{(K, [T])}a \Rightarrow (Eu)uB_{(K, [T])}b) \Rightarrow (Ez)zB_{(K, [T])}c$.

## 4. Conclusion

We have argued that Meta-theorem 2 is a valid Deduction Theorem of any first order theory. However, standard interpretations of Gödel's reasoning and conclusions are inconsistent with the consequences of this Meta-theorem in an arbitrary first order theory[16]. Hence, in the absence of constructive, and intuitionistically unobjectionable, reasons for denying the applicability of the Meta-theorem, and of its Corollary 2.1, to a first order theory, such interpretations ought not to be considered as definitive.

## Appendix 1: A model-theoretic proof of Corollary 2.1

We note that there is a model-theoretic proof of Corollary 2.1. The case $[T]|-_K [B]$ is straightforward.

If $[T]|-_K [B]$ does not hold, then, as noted in Meta-theorem 2, if $[T]|-_K [B]$ holds when we assume $[T]|-_K [A]$, then there is a sequence $<[B_1], [B_2], ..., [B_n]>$ of well-formed K-formulas such that $[B_1]$ is $[A]$, $[B_n]$ is $[B]$ and, for each $i > 1$, either $[B_i]$ is an axiom of K, or $[B_i]$ is in $[T]$, or $[B_i]$ is a direct consequence by some rules of inference of K of some of the preceding well-formed formulas in the sequence.

(We note that if $[T]$ is the set of well-formed K-formulas $\{[T_1], [T_2], ..., [T_l]\}$ then $([T] \& [A])$ denotes the well-formed K-formula $[T_1 \& T_2 \& ..., T_l \& A]$, and, $(T \& A)$ denotes its interpretation in M, i.e., $T_1 \& T_2 \& ..., T_l \& A$.)

If, now, any well-formed formula in $([T], [A])$ is false under an interpretation M of K, then $(T \& A) \Rightarrow B$ is vacuously true in M.

---

[16] See Appendix 2.

If, however, all the well-formed formulas in ([T], [A]) are true under interpretation in M, then the sequence <[B_1], [B_2], ..., [B_n]> interprets as a deduction in M, since the interpretation preserves the axioms and rules of inference of K (cf. [Me64], p57). Thus [B] is true in M, and so is (T & A) => B.

In other words, we cannot have ([T], [A]) true, and [B] false, under interpretation in M, as this would imply that there is some extension K' of K in which [T]|-$_{K'}$ [A], but not [T]|-$_{K'}$ [B]; this would contradict our hypothesis, which implies that, in any extension K' of K in which we have [T]|-$_{K'}$ [A], we also have [T]|-$_{K'}$ [B].

Hence, (T & A) => B is true in all models of K. By a consequence of Gödel's Completeness Theorem for an arbitrary first order theory ([Me64], p68, Corollary 2.15(*a*)), it follows that |-$_K$ ([T] & [A]) => [B]); and, ipso facto, that [T]|-$_K$ ([A] => [B]).

## Appendix 2: Gödel's reasoning and Corollary 2.1

In his seminal 1931 paper [Go31], Gödel meta-mathematically argues that, assuming any formal system of Peano Arithmetic, PA, is simply consistent, we can define an "undecidable" PA-proposition, [(A*x*)*R*(*x*)], such that (cf. [Go31], #1, p25):

   If [(A*x*)*R*(*x*)] is PA-provable, then [~(A*x*)*R*(*x*)] is PA-provable.[17]

Now, by Corollary 2.1, it should follow that:

   [(A*x*)*R*(*x*) => ~(A*x*)*R*(*x*)] is PA-provable,

and, therefore, that:

---

[17] Gödel essentially argues, number-theoretically, that, if the Gödel-number of [(A*x*)*R*(*x*)] is 17Gen*r*, and if this formula is PA-provable, then the PA-formula whose Gödel-number is Neg(17Gen*r*), i.e., [~(A*x*)*R*(*x*)], is also PA-provable if PA is assumed simply consistent.

[~(A*x*)*R*(*x*)] is PA-provable.

Since Gödel also proved that, if PA is assumed simply consistent, then [*R*(*n*)] is PA-provable for any, given, natural number *n*, Corollary 2.1 implies that PA is omega-inconsistent.

(We note that Gödel defined a first order theory K as *omega*-consistent if, and only if, for every well-formed formula [*F*(*x*)] of K, if |-$_K$[*F*(*n*)] for every numeral [*n*], then it is not the case that |-$_K$(E*x*)*F*(*x*) (cf. [Me64], p142; see also [Go31], p23-24).

However, this conclusion is inconsistent with standard interpretations of Gödel's reasoning, which, first, assert both [(A*x*)*R*(*x*)] and [~(A*x*)*R*(*x*)] as PA-unprovable, and, second, assume that PA can be omega-consistent[18]. Such interpretations, therefore, implicitly deny that the PA-provability of [~(A*x*)*R*(*x*)] can be inferred from the above meta-argument; ipso facto, they imply that Corollary 2.1 is false.

## References

[Ch37]   Church, A. 1937. Introduction to Mathematical Logic. Dover, New York.

[Go31]   Gödel, Kurt. 1931. *On formally undecidable propositions of Principia Mathematica and related systems I*. Translated by Elliott Mendelson. In M. Davis (ed.). 1965. The Undecidable. Raven Press, New York.

[Me64]   Mendelson, Elliott. 1964. Introduction to Mathematical Logic. Van Norstrand, Princeton.

[Sh67]   Shoenfield, J. R. 1967. Mathematical Logic. Association for Symbolic Logic, Urbana.

---

[18] We note that Gödel's Incompleteness Theorems assume significance only if we presume that the arithmetic, in which they are derived, *can* be *omega*-consistent (cf. [Go31], p23-24)